\documentclass[twocolumn, 10pt]{amsart}

\usepackage[english]{babel}
\usepackage{amsrefs}
\usepackage{hyperref}
\usepackage{amsmath,amssymb,amsfonts,latexsym, amsthm}
\usepackage[utf8]{inputenc}
\usepackage{xcolor}
\usepackage{empheq}
\usepackage[margin=2cm]{geometry}
\usepackage{mathtools}
\usepackage{graphicx}
\usepackage{units}
\usepackage{breqn}

\usepackage{mathptmx}
\usepackage{float}
\usepackage{indentfirst}
\usepackage{subfig}
\usepackage{tikz}
\usetikzlibrary{shapes,arrows,shadows}
\usetikzlibrary{decorations.markings}
\usepackage{color}
\usepackage{tikz-cd}

\graphicspath{ {images/} }
\usepackage[T1]{fontenc}
\usepackage{txfonts}
\usepackage{times}
\usepackage[all]{xy}
\usepackage{subfig}
\usepackage[normalem]{ulem}
\usepackage{hyperref}
\usepackage{wrapfig}
\usetikzlibrary{arrows}
\usepackage{pb-diagram}

\title{ A BRIEF TOUR THROUGH THE HISTORY OF
COMPLEX NUMBERS}
\author{Jhon Alexander Arredondo Garc\'ia and Camilo Ram\'irez Maluendas}

\address{Facultad de Matem\'aticas e Ingenier\'ias, Fundaci\'on Universitaria Konrad Lorenz, Bogot\'a, Colombia.}
\email{alexander.arrendo@konradlonrez.edu.co}
\address{Departamento de Matem\'atica y Estad\'istica\\ Facultad de Ciencias Exactas y Naturales\\ Universidad Nacional de Colombia, Sede Manizales \\
Manizales 170004, Colombia.}
\email{camramirezma@unal.edu.co}

\begin{document}

\maketitle

\section{Abstract}
In this paper, we chronologically recount several situations that have contributed to the development and formalization of the objects known as imaginary or complex numbers. We will begin by introducing the earliest documented knowing for calculating the square root of a negative quantity, attributed to the Greek mathematician Heron of Alexandria.  From there, we will progress through history to explore the formal concept of complex numbers given by William Rowan Hamilton.

\section{THE FIRST APPEARANCES}

The oldest known written evidence for calculating the square root of a negative quantity dates back to around 75 AD. This evidence can be found in the book  \emph{Stereometria}, written by the Greek mathematician, physicist, and engineer \emph{Heron of Alexandria}, who lived approximately between 10 and 85 AD. How did the square root of a negative quantity appear? Well, in such book the Greeks proposed to calculate the height of a  \emph{truncated pyramid with a square base and top} (see Figure \ref{Fig:Truncated pyramid}). In this context, $a=28$ and $b=4$ units represented the dimensions of the lower and upper squares, respectively, while $c=15$ units denoted the length of the slanted edge."

\begin{figure}[!ht]
		\begin{center}	
		\begin{tikzpicture}[baseline=(current bounding box.north)]
		\begin{scope}[scale=1.0]
			\clip (-2.2,-0.1) rectangle (2.2,4.2);
            \draw[fill=green!10, line width=1pt](0,0) -- (2,1) -- (0,2) -- (-2,1) -- (0,0);
            \draw[fill=yellow!10, line width=1pt](0,0) -- (2,1) -- (1,3.5) -- (0,3) -- (0,0);
           \draw[fill=yellow!10, line width=1pt](0,0) -- (0,3) -- (-1,3.5) -- (-2,1) -- (0,0);
            \draw[dashed, brown, line width=1pt](2,1) -- (0,2) -- (-2,1);
            \draw[dashed, brown, line width=1pt](0,4) -- (0,2);
           \draw[dashed, brown, fill=yellow!10, line width=0.8pt](0,0.5) -- (1,1) -- (0,1.5) -- (-1,1) -- (0,0.5);
           \draw[line width=1pt](0,0) -- (0,4);
            \draw[dashed, brown, line width=1pt](1,1) -- (2,1);
           \draw[dashed, brown, line width=1pt](1,1) -- (1,3.5);
           \draw[dashed, brown, line width=1pt](-1,1) -- (-2,1);
           \draw[dashed, brown, line width=1pt](-1,1) -- (-1,3.5);
          \draw[fill=yellow!10, line width=1pt](0,3) -- (1,3.5) -- (0,4) -- (-1,3.5) -- (0,3); 
           \draw[dashed, brown, line width=1pt](0,3) -- (0,4);
           \node[blue] at (1,0.2) {$a$};
           \node[blue] at (0.5,3) {$b$};
            \node[blue] at (-0.8,2) {$h$};
           \node[blue] at (1.3,1.15) {$x$};
          \node[blue] at (-1.7,2.3) {$c$};
			\end{scope}
			\end{tikzpicture} 
			 \caption{\emph{Truncated pyramid.}}
	\label{Fig:Truncated pyramid}
	\end{center}	
\end{figure}
We are familiar with the formula to calculate this height:
$$h=\sqrt{c^2-\left(\dfrac{a-b}{2}\right)^2}.$$ Heron was also aware of this formula. When he tackled this problem, he expressed it using modern notation:
\begin{eqnarray*}
h&=&\sqrt{\left(15\right)^2-\left(\dfrac{28-4}{2}\right)^2}=\sqrt{225-2(144)},\\
&=&\sqrt{225-144-144}=\sqrt{81-144}.\\
\end{eqnarray*}
From the previous equality one can conclude that the value of the height of the pyramid is $h=\sqrt{-63}$. However, in the subsequent step, Heron calculated $h=\sqrt{144-81}$ (refer to \cite{Nah}*{p. 4}, \cite{NRICH}, \cite{Smi}*{p. 261}). It's possible that a scribe made a writing error!

Another Greek mathematician who found himself in the same circumstances was the Hellenistic \emph{Diophantus of Alexandria}. Diophantus lived in Alexandria during the 3rd century AD. In his well-known treatise titled \emph{Arithmetica},which was written around 275 AD, Diophantus proposed a problem, which asked to find the measure of the sides of a right triangle with area $7$ square units and perimeter $12$ units (see Book VI, problem 22 \cite{Ras}). If we denote the lengths of the sides of this right triangle as $x$ and $y$, the problem can be expressed algebraically as follows:
\[
\frac{xy}{2}=7 \text{\quad and  \quad} x^2+y^2=(12-x-y)^2.
\] 
By applying the substitution method to the previous equations, we derive a second-degree equation
\begin{equation}\label{ec:cuadratica}
172x-24x^2-336=0.
\end{equation}
When we apply the general quadratic formula, we can determine that the solutions to equation \eqref{ec:cuadratica} are as follows:
\[
x=\frac{43\pm \sqrt{-167}}{12}.
\]
Nevertheless, Diophantus stated that equation \eqref{ec:cuadratica} could only be solved if half of the coefficient of 'x,' squared, minus the product of $24$ and $336$, results in a square. Clearly
\[
\left(\frac{172}{2}\right)^2-(24)(336)=-668,
\]
it is not a square! Because, it does not exist a real number multiplying by itself resulting in $-668$.

It seems that during Diophantus's time, the concept of an object being the square root of a negative quantity was not accepted. Around 850 BC, the Indian mathematician \emph{Mahaviracarya}, in his book \emph{Ganita-Sara-Sangraha}, established the following law to treat the roots of negative quantities: 

\begin{center}
	\begin{minipage}{0.9\linewidth}
		\vspace{5pt}
		{\small
		\emph{``As in the nature of things a negative (quantity) is not a square (quantity), it has  therefore no square root.''} \cite{Smi2}*{p. 312}.
		}
		\vspace{5pt}
	\end{minipage}
\end{center}

Likewise, the well-known Indian mathematician and astronomer \emph{Bhasker Acharij} (\textit{Bhaskara Acharya}, c. 1114-1185), mentioned in his algebra book titled \emph{Bijaganita} that: 
\begin{center}
	\begin{minipage}{0.9\linewidth}
		\vspace{5pt}
		{\small
		\emph{``But if any one asks the root of  negative $9$  I say the question is absurd, because never can be a negative square as has been shown.''}  \cite{Str}*{p. 15}.
		}
		\vspace{5pt}
	\end{minipage}
\end{center}

\section{first formalizations}

In the year 1545 the Italian polymath  \emph{Gerolamo Cardano} (1501-1576) published his work \emph{Ars magna}, in which by first time were introduced the algebraic object of the form  $a+\sqrt{-1}b$, where $a$ and $b$ are real numbers \cite{vander}*{p. 56}.  In his work, Cardano divulged  the solution to the cubic equations,  a method that had been developed in secret by the mathematicians \emph{Scipione del Ferro} (1465-1526) and \emph{Niccol\`o Fontana Tartaglia} (1501-1557). In addition,  he revealed the solution to the quartic equation, a method developed with \emph{Lodovoco Ferrari} (1522-1565). How do you solve a cubic polynomial? First, take the cubic equation
\[
ay^3+by^2+cy+d=0.
\]
Then substitute the new variable $x=y+\dfrac{b}{3a}$ into the above equation to obtain the reduced form
\begin{equation}\label{forma_reducida}
x^3+px+q=0,
\end{equation}
where
\[
p=\frac{3ac-b^2}{3a^2} \text{\quad and \quad} q=\frac{2b^3-9abc+27a^2d}{27a^3}.
\]
Now, replace $x=u+v$ in $x^{3}$ to obtain:
\begin{equation*}
\begin{split}
x^3& =(u+v)^3=u^3+v^3+3u^2v+3v^2u,\\
   &=u^3+v^3+3uv(u+v)=u^3+v^3+3uvx. 
\end{split}
\end{equation*}
This leads to the equation:
\begin{equation}\label{eq:comparar}
x^3-3uvx-u^3-v^3=0.
\end{equation}
Now, by comparing \eqref{forma_reducida} and \eqref{eq:comparar}, we obtain the following equations:
\begin{align}
3uv=-p,\label{eq:final_1}\\
u^3+v^3=-q.\label{eq:final_2}
\end{align}

In \eqref{eq:final_1}, express $v$ in terms of $p$ and $u$, then substitute this into \eqref{eq:final_2} to obtain:
\[
u^6+qu^3-\frac{p^3}{27}=0,
\] 
which is a quadratic equation in $z=u^{3}$. Rewrite the preceding equation as:
\[
z^2+qz-\dfrac{p^3}{27}=0.
\]
This leads to:
\[
z=u^3=-\frac{q}{2}\pm\sqrt{\frac{q^2}{4}+\frac{p^3}{27}}.
\]
Consider the positive root \footnote{Contrary, if you take the negative root, then you would obtain the same value to $x$.} to obtain:
\begin{equation*}\label{eq:u}
u=\sqrt[3]{-\frac{q}{2}+\sqrt{\frac{q^2}{4}+\frac{p^3}{27}}}.
\end{equation*}\label{eq:v}
Given that $v^3=-q-u^3$ (see equation \eqref{eq:final_2}), then
\begin{equation*}
v=\sqrt[3]{-\frac{q}{2}-\sqrt{\frac{q^2}{4}+\frac{p^3}{27}}}.
\end{equation*}
Finally, since $x=u+v$, it holds:
\begin{equation*}
x=\sqrt[3]{-\frac{q}{2}+\sqrt{\frac{q^2}{4}+\frac{p^3}{27}}}+\sqrt[3]{-\frac{q}{2}-\sqrt{\frac{q^2}{4}+\frac{p^3}{27}}}.
\end{equation*}
This is Cardano's method!

In Chapter 37 of the book \emph{Ars magna}, Cardano proposed the following problem: To find two numbers whose sum is $10$ and whose product is $40$. From algebraic standpoint, we need to determine the values $x$ and $y$ such that:
\[
x+y=10 \text{\, and \,} xy=40.
\]
After performing several basic operations, Cardano arrived at terms $5+\sqrt{-15}$ and $5-\sqrt{-15}$ as solutions \cite{Cardano}*{p. 219}. Then, he wrote:
\begin{center}
	\begin{minipage}{0.9\linewidth}
		\vspace{5pt}
		{\small
			\emph{``Putting aside the mental tortures involved, multiply $5+\sqrt{-15}$ and $5-\sqrt{-15}$, making $25 -(-15)$ which is $+15$. Hence this product is $40$... This is truly sophisticated.''} \cite{vander}*{p. 56}.
		}
		\vspace{5pt}
	\end{minipage}
\end{center}

The Italian mathematician  \emph{Rafael Bombelli} (1526-1572) also delved into cubic and quartic equations. In his book titled \emph{L' Algebra} \cite{Bombelli}, published in 1572, he referred to the object  $+\sqrt{-1}$ as \emph{pi\`u di meno} and the object $-\sqrt{-1}$ as \emph{meno di meno}. Moreover, he introduced eight rules for multiplying values involving the square negative root of $-1$ \cite{Fla}*{p. 25}. These rules are described in the Table \ref{Table:Bombelli}. It's worth noting that the symbol $\sqrt{-1}$ first emerged through the work of \emph{Albert Girard} (1595-1632) in his work \emph{Invention Nouvelle en l'Alg\`ebre}, published in 1629.
\begin{table*}[!ht]
\begin{center}
		\begin{tabular}{|l|r|}\hline
		\emph{1. Pi\`u  via pi\`u di meno  f\`a pi\`u di meno}.& $1 \cdot \sqrt{-1}=\sqrt{-1}$.\\
		\emph{2. Meno via pi\`u di meno f\`a meno di meno.} &$(-1)\cdot \sqrt{-1}=-\sqrt{-1}$.\\
		\emph{3. Pi\`u via  meno  di meno  f\`a meno di meno.}& $1\cdot\left(-\sqrt{-1}\right)=-\sqrt{-1}$.\\
		\emph{4. Meno via meno  di meno f\`a pi\`u di meno.} & $(-1)\cdot (-\sqrt{-1})=\sqrt{-1}.$\\
		\emph{5. Pi\`u di meno via pi\`u di meno f\`a meno.}& $\sqrt{-1}\cdot \sqrt{-1}=-1$. \\
		\emph{6. Pi\`u di meno via meno di meno f\`a pi\`u.}& $\sqrt{-1}\cdot (-\sqrt{-1})=\sqrt{-1}$. \\
		\emph{7. Meno di meno  via pi\`u di meno f\`a pi\`u.}& $\left(-\sqrt{-1}\right)\cdot \sqrt{-1}=1$.\\
		\emph{8. Meno di meno  via meno  di meno f\`a meno.} &$\left(-\sqrt{-1}\right)\cdot \left(-\sqrt{-1}\right)=-1$.\\ \hline
		\end{tabular}
	\caption{Rules of multiplication given by Bombelli.}
	\label{Table:Bombelli}
\end{center}
\end{table*}

In Chapter $2$, Bombelli found a solution to the equation $x^3=15x+4$ using Cardano's method. He obtained the root:
\[
x=\sqrt[3]{2+\sqrt{-121}}+\sqrt[3]{2-\sqrt{-121}},
\]
which he called \emph{sophisticated} root. Also, he recognized the solution $x=4$. Due to this, Bombelli became interested in the study of the cubic root of an imaginary number and aimed to found quantities $a$ and $b$ such that: 
\begin{equation}\label{eq:bom0}
\sqrt[3]{2+\sqrt{-121}}=a+\sqrt{-b}.
\end{equation}
How did Bombelli determine the values of $a$ and $b$?
From \eqref{eq:bom0}, he obtained:
\begin{equation*}
\begin{split}
2+\sqrt{-121}&=(a+\sqrt{-b})^3,\\
\end{split}
\end{equation*}
which expanded to:
\begin{equation*}
\begin{split}
&=a^3+(\sqrt{-b})^3+3a^2\sqrt{-b} +3a(\sqrt{-b})^2,\\
&=a^3-3ab + (3a^2-b)\sqrt{-b}.
\end{split}
\end{equation*}
By matching coefficients, he obtained the following equations:
\begin{align}
2&=a^3-3ab, \label{eq:bom1}\\
\sqrt{-121}&=(3a^2-b)\sqrt{-121}.\label{eq:bom2}
\end{align}
Then, he subtracted the two equations above and took cubic root to obtain:
\begin{equation}\label{eq:bom3}
\sqrt[3]{2-\sqrt{-121}}=a-b\sqrt{-1},
\end{equation}
He multiplied equations \eqref{eq:bom0} and \eqref{eq:bom3} to hold:
\begin{equation}\label{eq:bom4}
25=\sqrt[3]{125}=a^2+b.
\end{equation}
He then substituted \eqref{eq:bom4} into \eqref{eq:bom1} to obtain the cubic equation:
\[
4a^3-15a=2.
\]
This equation has the solution $a=2$. Finally, he substituted the value of $a$ into equation  \eqref{eq:bom1} to find $b=1$. Hence, Bombelli figured out the conjugated complex root\footnote{The term conjugate is due to Augustin-Louis Cauchy \cite{Green}*{p. 103}.} 
\[
\sqrt[3]{2+\sqrt{-121}}=2+\sqrt{-1} \text{\quad and \quad} \sqrt[3]{2-\sqrt{-121}}=2-\sqrt{-1},
\]
 whose sum is the value $4$.
 
 In 1637, the French philosopher, mathematician, and scientist  \emph{Ren\'e Descartes} (1596-1650) published the work titled \emph{La G\'eom\'etrie} \cite{Des}, which described his ideas about geometry applied to algebra. One of them involved associating lengths to line segments to find the length of a line segment from a quadratic equation. Specifically, in Book $I$ \cite{Des}*{p. 14}, Descartes aimed to construct a line segment whose length would satisfy the quadratic equation:
\[
x^2=ax-b^2,
\]
where $a$ and $b$ are positive real numbers. From the quadratic formula, we already  know the following solution
\[
x=\frac{1}{2}a\pm\sqrt{\frac{1}{4}a^2-b^2}.
\]
Nevertheless, Descartes resolved this puzzle by constructing two line segments, $\overline{NL}$ and $\overline{LM}$ with lengths $\frac{1}{2}a$ and $b$, respectively, as shown Figure \ref{segmentos}-a. He then drew the parallel straight line to $\overline{NL}$ passing through point $M$, and he took the circle with center $N$ and radius $NL=\frac{1}{2}a$. This resulted in the discovery that the circle and the parallel straight line intersected at two points, $Q$ and $R$, as shown in Figure \ref{segmentos}-b. Consequently, Descartes concluded that the length of the line segments $\overline{MQ}$ and $\overline{MR}$ were the solutions of his quadratic equation, \emph{i.e.}, 
\[
MQ=\frac{1}{2}a-\sqrt{\frac{1}{4}a^2-b^2} \text{\, and \,} MR=\frac{1}{2}a+\sqrt{\frac{1}{4}a^2-b^2}.
\]
\begin{figure*}[!ht]
     \centering
     \begin{tabular}{ccc}
      \begin{tikzpicture}[baseline=(current bounding box.north)]
		\begin{scope}[scale=1.0]
			\clip (-0.4,-0.4) rectangle (3.2,6.2);                     
         \draw[brown, line width=1pt](0,0) -- (0,3);
         \draw[brown, line width=1pt](0,0) -- (1,0);
           \node[brown] at (-0.2,-0.2) {$L$};
          \node[brown] at (1.2,-0.2) {$M$}; 
         \node[brown] at (-0.2,3) {$N$};
        \draw[line width=0.8pt](0,0.3) -- (0.3,0.3) -- (0.3,0);
			\end{scope}
			\end{tikzpicture} &    & 
		\begin{tikzpicture}[baseline=(current bounding box.north)]
		\begin{scope}[scale=1.0]
			\clip (-0.4,-0.4) rectangle (3.2,6.2);                     
         \draw[brown, line width=1pt](0,0) -- (0,3);
         \draw[brown, line width=1pt](0,0) -- (1,0);
        \draw [blue, line width=1pt] (0,0) arc(-90:90:3);
          \draw[red, line width=1pt](1,0) -- (1,5.84);
           \node[brown] at (-0.2,-0.2) {$L$};
          \node[brown] at (1.2,-0.2) {$M$}; 
         \node[red] at (1.2,0.5) {$Q$};
            \node[red] at (1.2,6) {$R$};
         \node[brown] at (-0.2,3) {$N$};
        \draw[line width=0.8pt](0,0.3) -- (0.3,0.3) -- (0.3,0);
			\end{scope}
			\end{tikzpicture} \\
\text{a. \emph{Line segments $\overline{NL}$ and $\overline{LM}$.}} & &\text{b. \emph{Common points  $Q$ and $R$.}}\\
     \end{tabular}
     \caption{\emph{Geometric construction of the root of the quadratic equation $x^2=ax-b^2$.}}
     \label{segmentos}
\end{figure*}
Why? Descartes used the \emph{Power point Theorem}; which states: If a tangent line segment and a secant line segment have a common point outside a circle, then the square of the tangent line segment is equal to the product of the secant line segment and the exterior line segment. Consequently, Descartes deduced: 
\begin{equation}\label{eq:potencia_punto}
MQ \cdot MQ = (LM)^2.
\end{equation}

Based on this last equality, he propose the following reasoning: if the length of the line segment $\overline{MQ}$ is $MQ=x$, then the length of the line segment $\overline{MR}$ must be $MR=a-x$. Therefore, the equation \eqref{eq:potencia_punto} is rewritten as: 
\[
x(a-x)=b^2 \text{\, or \,} x^2=ax-b^2.
\]
From these equations, Descartes deduced that the length of the segment line $\overline{MQ}$ also satisfies the equation  $x^2=ax-b^2$. 

Contrary, if the length of the line segment $\overline{MR}$ is $MR=x$, then the length  of the line segment $\overline{MQ}$ must be $MQ=a-x$. Then, the equation \eqref{eq:potencia_punto} is rewritten in the form: 
\[
x(a-x)=b^2 \text{\, or \,} x^2=ax-b^2.
\]
From the above equations, Descartes deduced that the length of the line segment $\overline{MR}$ is solution  of the equation $x^2=ax-b^2$. 

It's important to note that Descartes's construction relies on the inequality $\frac{1}{2}a\geq b$. If this inequality were reversed, the circle and the straight line passing through the points $R$ and $Q$ would not intersect. In this case, the roots of the equation would be of the form $a \pm \sqrt{-b}$. From this remark, Descartes wrote:

\begin{center}
	\begin{minipage}{0.9\linewidth}
		\vspace{5pt}
		{\small
			\emph{``Et si le cercle, qui ayant son centre au point $N$ passe par le point $L$, ne coupe ni ne touche la ligne droite $MQR$, il n'y a aucune racine en l'\'equation, de fa\c con qu'\ on peut assurer que la construction du probl\`eme propos\'e est impossible.''} \cite{Des}*{p. 15}.
		}
		\vspace{5pt}
	\end{minipage}
\end{center}
Due to this, the sophisticated roots began to be called \emph{impossible} numbers. However, Descartes coined the term \emph{imaginary} in this Book III:
\begin{center}
	\begin{minipage}{0.9\linewidth}
		\vspace{5pt}
		{\small
			\emph{``...on ne f\c caurait les rendre autres qu'\        imaginaires.''} \cite{Des}*{p. 174}.
		}
		\vspace{5pt}
	\end{minipage}
\end{center}

It might seem that Descartes tried to find a geometric interpretation for imaginary or impossible quantities, but his works suggest that he was not successful in this regard. Another attempt to geometrically construct the quantity $\sqrt{-1}$ was made by the English clergyman and mathematician \emph{John Wallis}\footnote{He also introduced the symbol known as infinity $\infty$.} (1616--1703). In Chapter LXVI of his work titled \emph{Treatise of Algebra}, published in 1685, knowing the negative numbers, Wallis associated to each point of a straight line a unique real number and vice versa. This geometric representation corresponds to the real number line as we know it today. It seems that Wallis had a geometric concept of where to construct imaginary numbers, as he wrote:
\begin{center}
	\begin{minipage}{0.9\linewidth}
		\vspace{5pt}
		{\small
			\emph{``Where $\sqrt{}$ implies a Mean Proportional between a Pofitive and a Negative Quantity. For like  as $\sqrt{bc}$ fignifies a Mean Proportional between $+b$ and $+c$; or between $-b$ and $-c$; (either of which, by Multiplication, makes $+bc$) tween  $-b$ and $+c$; either of which  being Multiplied, makes $-bc$. And this as to Algebraick confiderantion, is the true notion of fuch Imaginary Root, $\sqrt{-bc}$.''} \cite{Wallis}*{p. 290}.
		}
		\vspace{5pt}
	\end{minipage}
\end{center}

In Chapter LXVII, Wallis introduced the quantity $\sqrt{-81}$ as the length of the line segment $\overline{BC}$ close to the perpendicular axis $\overline{PC}$ (see Figure \ref{Wallis}), which we know as the imaginary axis.
 \begin{figure}[!ht]
	\begin{center}
		\begin{tikzpicture}[baseline=(current bounding box.north)]
		\begin{scope}[scale=1.0]
		\clip (-4.2,-0.5) rectangle (4.2,4.4);   
         \draw [black, line width=1.0pt](-4,0) -- (4,0);
         \draw [black, line width=1.0pt](0,4) -- (-4,0);
         \draw [red, line width=1.0pt](0,4) -- (-1.5,0.65);
         \draw [dashed, red!50, line width=1.0pt](-1.5,0) -- (-1.5,0.65);
         \draw [dashed, red!50, line width=1.0pt](0,0) -- (-1.5,0.65);
         \draw [red, line width=1.0pt](0,4) -- (1.5,0.65);
         \draw [dashed, red!50, line width=1.0pt](1.5,0) -- (1.5,0.65);
         \draw [dashed, red!50, line width=1.0pt](0,0) -- (1.5,0.65);
          \draw [dashed, red!50, line width=1.0pt](4,0) -- (1.5,0.65);
         \draw [dashed, red!50, line width=1.0pt](0,0) -- (0,4);
         \draw[dashed, blue, line width=1pt] (0,2) circle(2);
         \node at (1.5,-0.3) {$\beta$}; 
        \node at (0,-0.3) {$C$}; 
        \node at (-4,-0.3) {$A$}; 
        \node at (-1.5,-0.3) {$\beta$};
        \node at (0,4.2) {$P$};
        \node at (4,-0.3) {$\alpha$}; 
        \node at (-1.8,0.75) {$B$};
        \node at (1.8,0.75) {$B$};
			\end{scope}
			\end{tikzpicture}\\
		\caption{Wallis's construction.}
		\label{Wallis}
	\end{center}
\end{figure}
Wallis's work is a significant contribution to the geometric understanding of complex numbers.

\section{The intervention of the heavyweights}

One of the architects of the \emph{infinitesimal calculus} also studied imaginary roots. The German  mathematician, philosopher, scientist, and diplomat \emph{Gottfried Wilhelm Leibniz} (1646-1716) proved in 1676 the following equality
\[
\sqrt{1+\sqrt{-3}}+\sqrt{1-\sqrt{-3}}=\sqrt{6}.
\]
In addition, in 1702 (see \cite{Smi}*{p. 264}) he  expressed the term $x^4+a^4$ in the form 
\[
x^4+a^4=\left(x\pm a\sqrt{-\sqrt{-1}}\right)\left(x\pm a\sqrt{\sqrt{-1}}\right).
\]
By employing trigonometric substitution, we can prove that the Primitive Function of $\dfrac{1}{x^2+1}$ is  $\arctan (x)$. Using integral calculus notation:
\[
\int \frac{1}{x^2+1}dx=\arctan (x)+ C.
\]

In 1702, the Swiss mathematician \emph{Johann Bernoulli} (1667-1748) rewrote the rational function $\dfrac{1}{x^2+1}$ as the sum of other two ``simple'' rational functions 
\[
\frac{1}{x^2+1}=\frac{1}{2\sqrt{-1}} \left(\frac{1}{x-\sqrt{-1}}-\frac{1}{x+\sqrt{-1}}\right).
\]
Additionally, he introduced the \emph{Inverse Tangent Function} (up constant) as a \emph{Complex Logarithm Function} \footnote{From the historical point of view,  it can be inferred that this logarithm originated there.}
\begin{equation}\label{eq:bernoulli}
\begin{split}
\arctan(x) &=\int\frac{1}{x^2+1}dx,\\
&=\frac{1}{2\sqrt{-1}}\int \left(\frac{1}{x-\sqrt{-1}}-\frac{1}{x+\sqrt{-1}}\right)dx,\\ 
&=\frac{1}{2\sqrt{-1}}\ln \left( \frac{x-\sqrt{-1}}{x+\sqrt{-1}}\right).
\end{split}
\end{equation}
Twelve years later, in March 1714, the  English mathematician \emph{Roger Cotes} (1682-1716) published the manuscript titled  \emph{Logometria} in the \emph{Philosophical Transactions of Royal Society} \cite{Gow}. In this work, Cotes introduced  the following interesting formula
\begin{equation}\label{eq:Cotes}
\ln\left(\cos (\phi)+\sqrt{-1}\sin(\phi)\right )=\sqrt{-1}\phi,
\end{equation}
which is equivalent to the very well known Euler's formula  \[
\cos (\phi)+i\sin(\phi) =e^{i\phi}.
\]
How did Cotes obtain this formula? 

In modern notation, he considered the ellipse represented by the equation $\dfrac{x^2}{a^2}+\dfrac{y^2}{b^2}=1$, where $a$ and $b$ are the length of the semi axis. Then, he took the solid of revolution $S$, which is a piece of ellipsoid as shown the Figure \ref{elipse}-b. This solid is obtained by rotating the curve $\gamma$ around the $y$-axis, where $\gamma$ is the intersection of the ellipse and first quadrant of the Euclidean plane (see Figure \ref{elipse}-a). 
\begin{figure*}[!ht]
     \centering
     \begin{tabular}{ccc}
      \begin{tikzpicture}[baseline=(current bounding box.north)]
		\begin{scope}[scale=1.0]
			\clip (-1.7,-0.4) rectangle (3.7,2.2);   
          \draw [->, >=latex, black!30, line width=1.0pt](-1.6,0) -- (3.6,0); 
         \draw [->, >=latex, black!30, line width=1.0pt](0,-0.2) -- (0,2);                         
        \draw [blue, line width=1pt] (3,0) arc
		 		[
		 		start angle=0,
		 		end angle=90,
		 		x radius=30mm,
		 		y radius =15mm
		 		] ;
\draw [->, >=latex, line width=1pt] (0,2) arc
		 		[
		 		start angle=90,
		 		end angle=-250,
		 		x radius=4mm,
		 		y radius =1mm
		 		] ;
           
          \node[blue] at (3,-0.2) {$a$}; 
         \node[blue] at (-0.2,1.5) {$b$};
			\end{scope}
			\end{tikzpicture} &    &  \begin{tikzpicture}[baseline=(current bounding box.north)]
		\begin{scope}[scale=1.0]
			\clip (-3.7,-0.6) rectangle (3.7,2.2);   
         \draw [blue, line width=1.0pt](0,-0.5) -- (0,1.5);                         
        \draw [blue, line width=1pt] (3,0) arc
		 		[
		 		start angle=0,
		 		end angle=180,
		 		x radius=30mm,
		 		y radius =15mm
		 		] ;
    \draw [blue, line width=1pt] (3,0) arc
		 		[
		 		start angle=0,
		 		end angle=-180,
		 		x radius=30mm,
		 		y radius =5mm
		 		] ;           
    \draw [blue, dashed, line width=1pt] (3,0) arc
		 		[
		 		start angle=0,
		 		end angle=180,
		 		x radius=30mm,
		 		y radius =5mm
		 		] ;       
 \draw [blue, line width=1pt] (2.6,0.7) arc
		 		[
		 		start angle=0,
		 		end angle=-180,
		 		x radius=26mm,
		 		y radius =4mm
		 		] ; 
 \draw [blue, dashed, line width=1pt] (2.6,0.7) arc
		 		[
		 		start angle=0,
		 		end angle=180,
		 		x radius=26mm,
		 		y radius =4mm
		 		] ;   
 \draw [blue, line width=1pt] (2,1.1) arc
		 		[
		 		start angle=0,
		 		end angle=-180,
		 		x radius=20mm,
		 		y radius =2mm
		 		] ;   
 \draw [blue, dashed, line width=1pt] (2,1.1) arc
		 		[
		 		start angle=0,
		 		end angle=180,
		 		x radius=20mm,
		 		y radius =2mm
		 		] ; 
 \draw [blue, line width=1pt] (0,1.5) arc
		 		[
		 		start angle=90,
		 		end angle=180,
		 		x radius=8mm,
		 		y radius =20mm
		 		] ;  
 \draw [blue, line width=1pt] (0,1.5) arc
		 		[
		 		start angle=90,
		 		end angle=0,
		 		x radius=8mm,
		 		y radius =20mm
		 		] ;     
 \draw [blue, line width=1pt] (0,1.5) arc
		 		[
		 		start angle=90,
		 		end angle=160,
		 		x radius=22mm,
		 		y radius =28mm
		 		] ; 
 \draw [blue, line width=1pt] (0,1.5) arc
		 		[
		 		start angle=90,
		 		end angle=20,
		 		x radius=22mm,
		 		y radius =28mm
		 		] ;               
			\end{scope}
			\end{tikzpicture}  \\
\text{a. \emph{Curve $\gamma$.}} & &\text{b. \emph{Piece of an ellipsoid.}}\\
     \end{tabular}
     \caption{\emph{Solid of revolution $S$.}}
     \label{elipse}
\end{figure*}

Later, he found the area of the solid of revolution $S$ and obtained the following relations:
\begin{equation}\label{area_1}
A(S)=\pi a\left(a +\frac{b^2}{\sqrt{a^2-b^2}}\ln\left(\frac{a+\sqrt{a^2-b^2}}{b}\right) \right), 
\end{equation}
if it satisfies the inequality  $a>b$, or
\begin{equation}\label{area_2}
A(S)=\pi a\left(a+\frac{b^2}{\sqrt{b^2-a^2}}\sin^{-1}\left(\frac{\sqrt{b^2-a^2}}{b}\right) \right), 
\end{equation}
if it satisfies the inequality $a<b$. Hence, Cotes concluded that both relations \eqref{area_1} and \eqref{area_2} were the same for both inequalities: $a>b$ and $a<b$, \emph{i.e,}
\begin{equation}\label{eq:igualda_cotes}
\begin{split}
\pi a\left(a +\frac{b^2}{\sqrt{a^2-b^2}}\ln\left(\frac{a+\sqrt{a^2-b^2}}{b}\right) \right)=&\\ \pi a\left(a+\frac{b^2}{\sqrt{b^2-a^2}}\sin^{-1}\left(\frac{\sqrt{b^2-a^2}}{b}\right) \right).
\end{split}
\end{equation}
Then Cotes defined the value $\phi$ as   
\[
\sin (\phi)= \dfrac{\sqrt{b^2-a^2}}{b}=\dfrac{\sqrt{-1}\sqrt{a^2-b^2}}{b} \text{\, and \,} \cos (\phi)=\dfrac{a}{b},
\] 
and, we rewrote the equation \eqref{eq:igualda_cotes} in the form
\begin{equation}
\begin{split}
\pi a\left(a +\frac{b^2}{\sqrt{a^2-b^2}}\ln\left( \cos(\phi)-\sqrt{-1}\sin(\phi)\right) \right)=&\\ \pi a\left(a+\frac{b^2}{\sqrt{-1}\sqrt{b^2-a^2}}\sin^{-1}\left(\sin(\phi)\right) \right).
\end{split}
\end{equation}
So that, he obtained the equality
\[
-\sqrt{-1}\phi=\ln\left(\cos(\phi)-\sqrt{-1}\sin(\phi)\right),
\]
which is equivalent to \eqref{eq:Cotes}.

Cotes's formula motivated the French mathematician \emph{Abraham de Moivre} (1667-1754) to introduce in 1730 his well-known formula, which links complex numbers and trigonometry 
\[
\left(\cos (\phi)+\sqrt{-1}\sin(\phi)\right)^{n}= \cos (n\phi) + \sqrt{-1}\sin (n\phi).
\] 

It's worth noting that if we consider  $\phi=\frac{\pi}{2}$ and substitute in the equation given by Cotes, then we obtain
\[
\frac{\pi}{2}=\frac{\ln\left(\sqrt{-1}\right)}{\sqrt{-1}},
\]
which is a nice representation of $\pi$ from the complex logarithm. This representation was known by \emph{Leonhard Paul Euler} (1707-1783) and \emph{Johann Bernoulli}!  \cite{Fla}*{p. 82}. 
On October 18, 1740, Euler wrote a latter to Johann Bernoulli, which showed the functions 
\[
y(x)=2\cos(x) \quad \text{\, and \,} \quad y(x)=e^{x\sqrt{-1}}+e^{-x\sqrt{-1}}.
\]
as the solutions of the differential equation
\[
y''+y=0, \quad y(0)=2 \text{\quad and \quad} y'(0)=0.  
\] 
So, Euler introduced the equations \footnote{The work by Euler having these results appeared in 1748 \cite{Nah}*{p. 143}.} 
\[
\sin(\phi)=\frac{e^{\sqrt{-1}\phi}-e^{-\sqrt{-1}\phi}}{2\sqrt{-1}} \text{ and } \cos(\phi)=\frac{e^{\sqrt{-1}\phi}+e^{-\sqrt{-1}\phi}}{2}.
\]
Euler also introduced the notation 
\[
a+b\sqrt{-1}=e^{C}(\cos(\phi)+\sqrt{-1}\sin(\phi))=e^{C}e^{\sqrt{-1}(\phi \pm 2\lambda\pi)},
\]
and defined the complex logarithm function
\[
\ln (a+b\sqrt{-1})=C+\sqrt{-1}(\phi \pm 2\lambda\pi),
\]
where $\lambda$ is a positive integer \cite{Kline}*{p. 411}. In addition, the notation $i$, which we currently use to designate the square root of -$1$, first appeared in the Euler's work titled \emph{De formulis differentialibus angularibus maxime irrationalibus, quas tamen per logarithmos et arcus circulares integrare licet. M. S. Academiae exhibit. die 5. Maii 1777} by Euler. In this work, he wrote:
\begin{center}
	\begin{minipage}{0.9\linewidth}
		\vspace{5pt}
		{\small
			\emph{``...formulam $\sqrt{-1}$ littera $i$  in posterum designabo,...''} 
		}
		\vspace{5pt}
	\end{minipage}
\end{center}

The problem of geometrically representing imaginary roots was tackled most successfully on March 10, 1797. At that time, the Danish–Norwegian mathematician \emph{Caspar Wessel} (1745-1818) submitted his work titled \emph{On the Analytical Representation of Direction. An Attempt Applied Chiefly to Solving Plane  and Spherical Polygons} \cite{Wessel} to the \emph{Royal Danish Academy of Sciences and Letters}. It was published two year later. In his work, Wessel geometrically represented each imaginary root $a+bi$ as \emph{directed line segment or vector}. He defined the algebra of vector, including addition and multiplication\footnote{Eventually, this work led to the name of the imaginary numbers.}.  From the geometric point of view, the addition is known as the \emph{Parallelogram Rule}. In the Figure \ref{vectores}-a is shown  the addition of vectors $\overrightarrow{AB}$ and $\overrightarrow{BC}$, which is denoted by $\overrightarrow{AB}$+$\overrightarrow{BC}$=$\overrightarrow{AC}$. 
\begin{figure*}[!ht]
	\centering
	\begin{tabular}{ccc}
		\begin{tikzpicture}[baseline=(current bounding box.north)]
		\begin{scope}[scale=1.0]
		\clip (-0.7,-0.6) rectangle (3.4,4);   
          \draw [->, >=latex, blue, line width=1.0pt](0,0) -- (2,1);                                 
           \draw [->, >=latex, green, line width=1.0pt](0,0) -- (3,3.5);   
          \draw [->, >=latex, red, line width=1.0pt](2,1) -- (3,3.5);     
          \node at (0,-0.2) {$A$}; 
         \node at (3,3.7) {$B$};
          \node at (2.3,1) {$C$};
			\end{scope}
			\end{tikzpicture} &    &  \begin{tikzpicture}[baseline=(current bounding box.north)]
		\begin{scope}[scale=1.0]
		\clip (-0.7,-0.5) rectangle (3.4,5.2);   
          \draw [->, >=latex, black!30, line width=1.0pt](-0.6,0) -- (3.2,0);
          \draw [->, >=latex, blue, line width=1.0pt](0,0) -- (2,1); 
          \draw[blue] (1,0) arc (0:27:1);                                
           \draw [->, >=latex, red, line width=1.0pt](0,0) -- (2,2);
           \draw[red] (1.5,0) arc (0:45:1.5);        
          \draw [->, >=latex, green, line width=1.0pt](0,0) -- (1.6,5);    
        \draw[green] (2,0) arc (0:72:2);   
          \node at (0,-0.3) {$\textbf{0}$}; 
         \node at (3,-0.3) {$E$}; 
         \node[red] at (2.3,2) {$B$};
          \node[blue] at (2.3,1) {$A$};
        \node[green] at (1.9,5) {$C$};
			\end{scope}
			\end{tikzpicture} \\
		\text{a. \emph{The Parallel Rule.}} & &\text{b. \emph{Multiplication of vectors.}}\\
	\end{tabular}
	\caption{\emph{Algebra of vectors.}}
	\label{vectores}
\end{figure*}
To define the multiplication of vectors, Wessel used a straight line (axis) as a guide. He drew two point on it, the origin  $\textbf{0}$ and $E$. Then, he defined the multiplication of the vectors $\overrightarrow{\textbf{0}A}$ and $\overrightarrow{\textbf{0}B}$ as the new vector $\overrightarrow{\textbf{0}C}$, which has a  length equal to the product of  product of the lengths of the vectors $\overrightarrow{\textbf{0}A}$ and $\overrightarrow{\textbf{0}B}$. Moreover, the direction of $\overrightarrow{\textbf{0}C}$ is the sum of the directions of the vectors $\overrightarrow{\textbf{0}A}$ and $\overrightarrow{\textbf{0}B}$. This means that the vector $\overrightarrow{\textbf{0}C}$ must have a length of $\left| \overrightarrow{\textbf{0}A}\right| \cdot \left|  \overrightarrow{\textbf{0}B}\right|$\footnote{The symbol  $|\, |$ denotes length or norm of a vector.} and a direction of $m(\angle C\textbf{0}E)=m(\angle A\textbf{0}E) +m(\angle B\textbf{0}E)$\footnote{The symbol $m$ denotes direction.}. See Figure \ref{vectores}-b. 

How did Wessel discover the multiplication of complex numbers? He did so based on the following remarks described in Table \ref{Table:1}.
 \begin{table*}[!ht]
	\begin{center}
		\begin{tabular}{|l|l|}\hline
			Multiplication& Geometry \\ \hline\hline
			$(1)\cdot (1)=1$&   Length $1$ and direction $0^{\rm o}.$\\ \hline
			$(1)\cdot (-1)=-1$&  Length $1$ and direction $180^{\rm o}$.\\ \hline
			$(-1)\cdot (-1)=-1$& Length $1$ and direction $0^{\rm o}$.   \\ \hline
			$(1)\cdot (\sqrt{-1})=\sqrt{-1}$& Length $1$ and direction $90^{\rm o}$.\\ \hline
			$(1)\cdot (-\sqrt{-1})=-\sqrt{-1}$& Length $1$ and direction $0^{\rm o}+270^{\rm o}=270^{\rm o}$.\\ \hline
			$(-1)\cdot (\sqrt{-1})=-\sqrt{-1}$& Length $1$ and direction $180^{\rm o}+90^{\rm o}=270^{\rm o}$.\\ \hline
			$(-1)\cdot (-\sqrt{-1})=\sqrt{-1}$& Length $1$ and direction $180^{\rm o}+270^{\rm o}=360^{\rm o}+90^{\rm o}=90^{\rm o}$.\\ \hline
			$(\sqrt{-1})\cdot (\sqrt{-1})=-1$ & Length $1$ and direction $90^{\rm o}+90^{\rm o}=180^{\rm o}$.\\ \hline
			$(\sqrt{-1})\cdot -(\sqrt{-1})=1$& Length $1$ and direction $90^{\rm o}+270^{o}=360^{\rm o}=0^{\rm o}$.\\ \hline
			$(-\sqrt{-1})\cdot (\sqrt{-1})=1$ & Length $1$ and direction $270^{\rm o}+90^{\rm o}=360^{\rm o}=0^{\rm o}.$\\ \hline
			$(-\sqrt{-1})\cdot -(\sqrt{-1})=-1$ & Length $1$ and direction $270^{\rm o}+270^{\rm o}=360^{\rm o}+180^{\rm o}=180^{\rm o}.$ \\ \hline
		\end{tabular}
	\caption{Wessel's remark.}
	\label{Table:1}
	\end{center}
	\end{table*}
How did Wessel associate each imaginary number $a+\sqrt{-1}b$ with a vector?  He accomplished this by drawing two straight lines on the Euclidean plane, which intersected at a point $\textbf{0}$ called \emph{origin}. These lines formed four right angles, representing the real and imaginary axes\footnote{Several historians attribute to Wessel the introduction of the perpendicular axis, known as the imaginary axis, alongside the real axis \cite{Nah}*{p. 53}.}. Then, he drew on such axis the points labeled as $1$, $-1$, $\sqrt{-1}$ and $-\sqrt{-1}$ (see Figure \ref{wessel}). 
 \begin{figure}[!ht]
	\begin{center}
		\begin{tikzpicture}[baseline=(current bounding box.north)]
		\begin{scope}[scale=1.0]
		\clip (-1.8,-1.8) rectangle (3.1,3.1);   
         \draw [->, >=latex, black!30, line width=1.0pt](-1.7,0) -- (3,0);
         \draw [->, >=latex, black!30, line width=1.0pt](0,-1.7) -- (0,3);
         \draw [->, >=latex, blue, line width=1.0pt](0,0) -- (2.5,2.5);
         \draw [blue, dashed, line width=1.0pt](0,2.5) -- (2.5,2.5)  -- (2.5,0);
         \node[red] at (0,1.5) {{\tiny$\bullet$}}; 
        \node[red] at (0,-1.5) {{\tiny$\bullet$}}; 
        \node[red] at (1.5,0) {{\tiny$\bullet$}}; 
        \node[red] at (-1.5,0) {{\tiny$\bullet$}}; 
        \node[blue] at (2.7,2.5) {$A$}; 
       \node[blue, rotate=45] at (1.2,1.6) {$a+\sqrt{-1}b$}; 
         \node at (1.5,-0.3) {$1$}; 
        \node at (-0.3,-0.3) {$0$}; 
        \node at (-1.5,-0.3) {$-1$}; 
         \node at (-0.8,-1.5) {$-\sqrt{-1}$};
         \node at (-0.6, 1.5) {$\sqrt{-1}$};
        \node[blue] at (2.5,-0.3) {$a$}; 
        \node[blue] at (2.5,0) {${\tiny \bullet}$}; 
        \node[blue] at (-0.3,2.5) {$b$}; 
        \node[blue] at (0,2.5) {${\tiny \bullet}$}; 
			\end{scope}
			\end{tikzpicture}\\
		\caption{\emph{Wessel's construction of $a+\sqrt{-1}b$.}}
		\label{wessel}
	\end{center}
\end{figure}
As a result, the imaginary $a+\sqrt{-1}b$ was identified as the vector $\overrightarrow{\textbf{O}A}$  (see Figure \ref{wessel}). This vector has length  $r:=\sqrt{a^2+b^2}$ and direction
\[
\phi= \left\{
\begin{array}{ ll }
\arctan \left(\dfrac{y}{x}\right)-\pi; & \text{ if } y<0,\, x<0; \\
-\dfrac{\pi}{2}; & \text{ if } y<0, \, x=0;\\
\arctan\left(\dfrac{y}{x}\right); & \text{ if } x>0; \\
\dfrac{\pi}{2}; & \text{ if } y>0, \, x=0;\\
\arctan\left(\dfrac{y}{x}\right)+\pi; & \text{ if } y\geq 0,\, x<0. \\
\end{array}
\right.
\]
From his geometric ideas, Wessel wrote the polar system
\begin{align*}
 a&=r\cos(\phi),\\
 b&=r\sin(\phi),\\
 a+\sqrt{-1}b&=r\left(\cos(\phi)+\sqrt{-1}\sin(\phi)\right).
\end{align*}
From this, he remarked that if you have the relation: 
\[
c+\sqrt{-1}d=r'\left(\cos(\alpha)+\sqrt{-1}\sin(\alpha)\right),
\]
then the multiplications of two imaginary numbers must be given by   
\begin{equation*}
\begin{split}
\left(a+\sqrt{-1}b\right)\cdot \left(c+\sqrt{-1}d\right)  &=rr'(\cos(\phi+\alpha)\\ &+\sqrt{-1}\sin(\phi+\alpha)),\\ 
&=\left(ac-bd +\sqrt{-1}(ad+bc)\right).
\end{split}
\end{equation*}
And, the addition of two imaginary numbers must be represented by
\begin{equation*}
\begin{split}
\left(a+\sqrt{-1}b\right)+ \left(c+\sqrt{-1}d\right)  &=r\cos(\phi)+ r'\cos(\alpha)\\ &+ 
\sqrt{-1}\left(r\sin{\phi}+r'\sin(\alpha)\right),\\ 
&=a+c +\sqrt{-1}\left(b+d\right).
\end{split}
\end{equation*} 

Who was \emph{Jean-Robert Argand}? He was an accountant in Paris, born in 1768 and passing away in 1822.  Unfortunately, very little is known about his life. Nevertheless, what we do know is that Argand wrote a pamphlet in 1806, titled \emph{Essai sur une mani\`{e}re de pr\'esenter les quantit\'es imaginaires dans les constructions g\'eom\'etriques} \cite{Argand}, although the authorship was initially unknown. This manuscript was read by the renowned French mathematician \emph{Adrien-Marie Legendre} (1752--1833). Legendre, in turn, mentioned it through a letter to \emph{Francois Francais} (1768-1810), who was a professor of mathematics with extensive military background. Upon Francais's death, his brother \emph{Jaques} (1775-1833) inherited his papers. Jaques also had extensive military background and mathematical knowledge. He was a professor at the \emph{Ecole Imp\'eriale d'\ Application du G\'enie et de l'\ Artillerie} in Metz. Indeed, Jaques found Legendre's letter addressed to his brother in which he described Argand's mathematical results. In this missive Legendre does not mention Argand. Motivated by the ideas presented in that letter, Jaques published an article in 1813 in the journal \emph {Annales de Math \'emathiques}, in which he introduced the fundamental concepts of complex geometry. In the last paragraph of his article, Jacques acknowledged the letter from Legendre to his deceased brother and expressed a desire to identify the anonymous author. It was at this point that Argand became aware of the situation and communicated with Jacques Francais. Then, in the next issue of the magazine, Francais reported that Argand was the first person to have developed the geometry of imaginary numbers. Neither of them had heard of Wessel's work on these subjects \cite{Nah}*{p. 74}. 

In his article, Argand provided an interpretation of the imaginary numbers, such as $a+b\sqrt{-1}$, as a point in the plane. He also described the geometric rules of multiplication and addition of imaginary numbers. In addition, Argand interpreted $\sqrt{-1}$ as a rotation by $180^{\rm o}$.

The prince of mathematics \emph{Karl Friedrich Gauss} (1777-1855) also made valuable contributions to the field, particularly concerning what were then referred to as imaginary numbers. In April 1831, Gauss renamed it as \emph{complex numbers} \cite{Gauss1}*{p. 102}. In 1799, he presented his doctoral thesis titled \emph{Demonstratio nova theorematis omnem functionem algebraicam rationalem integram unius variabilis in factores reales primi vel secundi gradus resolvi posse} \cite{Gauss}, in which he introduced the first of four proofs, published during his lifetime, of the Fundamental Theorem of Algebra. What does this theorem tell us?

The \textbf{Fundamental Theorem of Algebra} states that \emph{any polynomial function $p(z)=a_{n}z^{n} +\ldots+a_{0}$ of degree $n\geq 1$ has a root in the complex plane $\mathbb{C}$.} 

In 1816, Gauss introduced two new proofs for this theorem in his paper titled \emph{Demonstratio nova altera theorematis omnem functionem algebraicam
	rationalem integram unius variabilis in factores reales primi vel secundi gradus resolvi posse}.
Comm. Recentiores (Gottingae), \textbf{3}:107-142, 1816. In Werke III, 31-56. Marsden, J. E. in the book \emph{Basic Complex Analysis} p. 151  presented the proof of the Fundamental Theorem of Algebra using essentially the ideas provided by Gauss in this publication. Basically, the strategy involved proving the assertion by contradiction. Gauss assumed that the polynomial $p(z)$ had no roots and then constructed the new function $f(z)=\frac{1}{p(z)}$, the which turned out to be whole and bounded. Then, he applied Liouville's Theorem and deduced that the function $f$ must be constant \emph{ergo}, the polynomial $p(z)$ was also a constant. This fact was a contradiction, since the polynomial $p(z)$ was not constant. This reasoning led him to conclude that there must exist an element $z_0\in\mathbb{C}$ such that $p(z_0)=0$.

In April 1831, Gauss presented \cite{Gauss1}*{p. 102} his geometric ideas on complexes numbers to the Royal Society of Göttingen, coinciding with those of Argand. It seems that Gauss had already developed these notions as early as 1796, which predates Wessel's work. However, like many of his works, he did not publish them until his ideas had been matured \cite{Nah}*{p. 82}.

With the French mathematician, engineer, and physicist \emph{Augustin-Louis Cauchy} (1789--1857) began the study of functions with complex values, that is, functions with correspondence rule $f(x+ iy)$. In 1814, he presented his work titled \emph{M\'emoire sur les integrales d\'efinies} to the \emph{French Academy of Sciences}, which  contained his contributions to development of the theory of complex functions \cite{Ett}. In this memoirs, Cauchy proved that if a complex function $f$ is analytic over a region $G$ ``with suitable characteristics'', then the integral of $ f$ along the boundary of any rectangle contained in $G$ is zero. Cauchy also studied complex functions that are discontinuous at isolated points and obtained a formula, which is the essence of the calculus of residues.

Who introduced the definition of complex numbers that we usually learn in an elementary complex analysis course? Such definition first appeared in the work titled \emph{Theory of Conjugate Functions
	or Algebraic Couples: with a Preliminary Essay on Algebra as a Science of
	Pure Time} \cite{Ham} by \emph{William Rowan Hamilton} (1805--1865), published in  1837. Hamilton introduced ordered pairs  $(a,b)$, where $a$ and $b$ are real numbers, and defined the addition and the multiplication of these ordered pairs as follows
\[
\begin{array}{ccl}
(a,b)+(c,d)& := & (a+c,b+d);\\
(a,b)\cdot (c,d) & := & (ac-bd,bc+ad).
\end{array}
\]
Moreover, Hamilton defined the square root for any ordered pair $(a,b)\in\mathbb{R}^2$. He then identified  each real number $a$ as the ordered pair $(a,0)$ and computed 
\[
\sqrt{-1}=\sqrt{(-1,0)}=(0,1).
\]
From the preceding equality, he wrote 
\begin{center}
	\begin{minipage}{0.9\linewidth}
		\vspace{5pt}
		{\small
			\emph{``...be concisely denoted as follows, $\sqrt{-1}=(0, 1)$. In the theory of single numbers, the symbol $\sqrt{-1}$ is absurd, and denotes an impossible extraction, or a merely imaginary number; but in the theory of couples, the same symbol $\sqrt{-1}$ is significant, and denotes a possible extraction, or a real couple, namely (as we have just now seen) the principal square-root of the couple $(-1,0)$.''}
		}
		\vspace{5pt}
	\end{minipage}
\end{center}

What has happened after Hamilton? There are still too many fascinating anecdotes left to share. However, we invite readers to tell their \emph{imaginary tale}.
 
\vspace{4mm}
 
\textbf{ACKNOWLEDGEMENTS.} Camilo Ram\'irez Maluendas expresses his gratitude to Universidad Nacional de Colombia, Sede Manizales. He has dedicated this work to his beautiful family: Marbella and Emilio, in appreciation of their love and support.


\end{document}